\documentclass{gtart_h}

\def\ifplaintex{\expandafter\ifx\csname documentclass\endcsname\relax}

\def\gtp{{\mathsurround=0pt\it $\cal G\mskip-2mu$eometry \&\ 
$\cal T\!\!$opology $\cal P\!$ublications}}  

\def\recd{{\small Received:\qua\receiveddate\ifx\reviseddate\relax
\else\qquad Revised:\qua\reviseddate\fi\par}} 

\def\lognumber#1{\def\thelognumber{#1}}
\def\volumenumber#1{\def\thevolumenumber{#1}}
\def\volumeyear#1{\def\thevolumeyear{#1}}
\def\papernumber#1{\def\thepapernumber{#1}}
\def\pagenumbers#1#2{\def\startpage{#1}\def\finishpage{#2}}
\def\published#1{\def\publishdate{#1}}

\def\received#1{\def\receiveddate{#1}}
\def\revised#1{\def\reviseddate{#1}}
\def\accepted#1{\def\accepteddate{#1}}

\def\asciiaddress#1{\def\theasciiaddress{#1}}
\def\asciiemail#1{\def\theasciiemail{#1}}

\long\def\asciiabstract#1{\long\def\theasciiabstract{#1}}

\let\\\par\let\thelognumber\relax\let\thevolumenumber\relax
\let\thepapernumber\relax\let\thevolumeyear\relax\let\startpage\relax
\let\finishpage\relax\let\publishdate\relax\let\receiveddate\relax
\let\reviseddate\relax\let\accepteddate\relax\let\theasciititle\relax
\let\theasciiauthors\relax\let\theasciiaddress\relax
\let\theasciiabstract\relax

\let\theasciiemail\relax

\ifplaintex
\font\logobig=cmssbx10 scaled 3836
\font\logomed=cmssbx10 scaled 2557
\else
\font\logobig=cmssbx10 scaled 4200
\font\logomed=cmssbx10 scaled 2800
\fi

\long\def\makeagttitle{   
\count0=\startpage
\agt\hfill      
\hbox to 45truept{\vbox to 0pt{\vglue -13truept{\logomed A\kern -.37em{\logobig 
T}\kern -.38em G}\vss}\hss}
\break
{\small Volume \thevolumenumber\ (\thevolumeyear)
\startpage--\finishpage\nl
Published: \publishdate}

\vglue .25truein

{\parskip=0pt\leftskip 0pt plus
1fil\def\\{\par\smallskip}{\Large\bf\thetitle}\par\medskip} \vglue
0.05truein

{\parskip=0pt\leftskip 0pt plus 1fil\def\\{\par}{\sc\theauthors}
\par\medskip}%
 
\vglue 0.03truein 

{\small\leftskip 25truept\rightskip 25truept{\bf Abstract}\stdspace\theabstract

{\bf AMS Classification}\stdspace\theprimaryclass
\ifx\thesecondaryclass\relax\else; \thesecondaryclass\fi\par
{\bf Keywords}\stdspace \thekeywords\par}\vglue 7truept

}   

\ifplaintex
\font\phead=cmsl9 scaled 950
\font\pnum=cmbx10 scaled 913
\font\pfoot=cmsl9 scaled 950
\headline{\vbox to 0pt{\vskip -4.5mm\line{\small\phead\ifnum
\count0=\startpage ISSN 1472-2739 (on-line) 1472-2747 (printed)
\hfill {\pnum\folio}\else\ifodd\count0\def\\{ }%
\ifx\theshorttitle\relax\thetitle\else\theshorttitle\fi\hfill{\pnum\folio}
\else\def\\{ and }{\pnum\folio}\hfill\ifx\theshortauthors\relax\theauthors
\else\theshortauthors\fi\fi\fi}\vss}}
\footline{\vbox to 0pt{\vglue 0mm\line{\small\pfoot\ifnum\count0=\startpage
\copyright\ \gtp\hfill\else
\agt, Volume \thevolumenumber\ (\thevolumeyear)\hfill\fi}\vss}}
\else
\font=cmsl9 scaled 1050
\font\lnum=cmbx10 
\font=cmsl9 scaled 1050
\makeatletter
\def\@oddhead{{\small\lhead\ifnum\count0=\startpage ISSN 1472-2739 
(on-line) 1472-2747 (printed)\hfill {\lnum\number\count0}\else\ifodd\count0
\def\\{ }\ifx\theshorttitle\relax \thetitle \else\theshorttitle\fi\hfill
{\lnum\number\count0}\else\def\\{ and }{\lnum\number\count0}
\hfill\ifx\theshortauthors\relax 
\theauthors\else\theshortauthors\fi\fi\fi}}\def\@evenhead{\@oddhead}
\def\@oddfoot{\small\lfoot\ifnum\count0=\startpage\copyright\ \gtp\hfill\else
\agt, Volume \thevolumenumber\ (\thevolumeyear)\hfill\fi}
\def\@evenfoot{\@oddfoot}
\makeatother
\fi
\let\maketitlepage\makeagttitle
\let\makeshorttitle\maketitlepage
\let\maketitle\maketitlepage

\newwrite\gtoutfile
\long\gdef\makeheadfile{  
{\def\\{, }\def\s{ }
\immediate\openout\gtoutfile head.xxx
\immediate\write\gtoutfile{Proxy-for: \ifx\theasciiauthors\relax
\theauthors\else\theasciiauthors\fi\s<\ifx\theasciiemail\relax\theemail\else\theasciiemail\fi>}
\immediate\write\gtoutfile{\noexpand\\}
\immediate\write\gtoutfile{Authors: \ifx\theasciiauthors\relax
\theauthors\else\theasciiauthors\fi}
{\def\\{ }\immediate\write\gtoutfile{Title: \ifx\theasciititle\relax
\thetitle\else\theasciititle\fi}}
\immediate\write\gtoutfile{Subj-class: GT or SG, GR etc}
\immediate\write\gtoutfile{MSC-class: \theprimaryclass\ifx\thesecondaryclass\relax\else, \thesecondaryclass\fi}
\immediate\write\gtoutfile{Journal-ref: Algebr. Geom. Topol. \thevolumenumber\s
(\thevolumeyear) \startpage-\finishpage}
\immediate\write\gtoutfile{Comments: Published by Algebraic and
Geometric Topology at}
\immediate\write\gtoutfile{\s\s\s  http://www.maths.warwick.ac.uk/agt/AGTVol\thevolumenumber/agt-\thevolumenumber-\thepapernumber.abs.html}
\immediate\write\gtoutfile{\noexpand\\}
\immediate\write\gtoutfile{}
\ifx\theasciiabstract\relax
\immediate\write\gtoutfile{\theabstract}\else
\immediate\write\gtoutfile{\theasciiabstract}\fi
\immediate\write\gtoutfile{}
\immediate\write\gtoutfile{\noexpand\\}
\immediate\write\gtoutfile{}
\immediate\closeout\gtoutfile}}  

\def\maketitlepage{\makeagttitle\makeheadfile}
\let\makeshorttitle\maketitlepage
\let\maketitle\maketitlepage

\lognumber{61}
\volumenumber{5}
\volumeyear{2005}
\papernumber{61}
\pagenumbers{1535}{1553}
\received{9 March 2005} 
\revised{20 October 2005}
\accepted{24 October 2005}
\published{16 November 2005}

\usepackage{amssymb,amsmath,amscd,graphicx,psfrag}
\usepackage[all,dvips]{xy}

\graphicspath{{Figures/}}

\def\psfraga <#1,#2> #3#4{%
\psfrag {#3}{\smash{\rlap{\kern #1 \raise #2\hbox{#4}}}}}
\psfraga <-1.5pt,0pt> {d}{{\normalsize$\scriptstyle d$}}
\psfraga <-2.5pt,0pt> {-d}{{\small$\scriptscriptstyle {-}d$}}
\psfraga <-2.5pt,0pt> {+d}{{\small$\scriptscriptstyle {+}d$}}
\psfraga <-1.5pt,0pt> {v}{{\normalsize$\scriptstyle v$}}
\psfraga <-1.5pt,0pt> {v'}{{\normalsize$\scriptstyle v'$}}
\psfraga <-1.5pt,-1pt> {v''}{{\normalsize$\scriptstyle v''$}}
\psfraga <-1.5pt,0pt> {h}{{\normalsize$\scriptstyle h$}}
\psfraga <-1pt,0pt> {a}{{\normalsize$\scriptscriptstyle \alpha$}}
\psfraga <-1pt,0pt> {b}{{\normalsize$\scriptscriptstyle \beta$}}
\psfraga <-1pt,0pt> {g}{{\normalsize$\scriptscriptstyle \gamma$}}

\def\figref#1{\hyperlink{#1anchor}{Figure~\ref*{#1}}}
\def\anchor#1{\noindent\hypertarget{#1anchor}{\smash{$\phantom{99}$}}}

\newtheorem{theo}{Theorem}[section]
\newtheorem{prop}[theo]{Proposition}
\newtheorem{lemme}[theo]{Lemma}

\theoremstyle{definition}
\newtheorem{defi}[theo]{Definition}

\newtheorem*{remarque}{Remark}
\newtheorem*{question}{Question}

\def\Q{{\mathbb Q}}
\def\R{{\mathbb R}}
\def\Z{{\mathbb Z}}

\def\A{{\mathcal A}}
\def\C{{\mathcal C}}
\def\H{{\mathcal H}}

\begin{document}


\title{A Jones polynomial for braid-like isotopies\\of oriented links and its categorification}
\shorttitle{A Jones polynomial for braid-like isotopies of links}

\authors{Benjamin Audoux\\Thomas Fiedler}
\address{Laboratoire E.Picard, Universit\'e Paul Sabatier\\Toulouse, France}
\asciiaddress{Laboratoire E.Picard, Universite Paul Sabatier\\Toulouse, France}
\asciiemail{audoux@picard.ups-tlse.fr, fiedler@picard.ups-tlse.fr}
\gtemail{\mailto{audoux@picard.ups-tlse.fr}{\rm\qua 
and\qua}\mailto{fiedler@picard.ups-tlse.fr}}

\begin{abstract}
A braid-like isotopy for links in $3$--space is an isotopy which uses
only those Reidemeister moves which occur in isotopies of braids. We
define a refined Jones polynomial and its corresponding Khovanov
homology which are, in general, only invariant under braid-like
isotopies.
\end{abstract}
\asciiabstract{%
A braid-like isotopy for links in 3-space is an isotopy which uses
only those Reidemeister moves which occur in isotopies of braids. We
define a refined Jones polynomial and its corresponding Khovanov
homology which are, in general, only invariant under braid-like
isotopies.}



\maketitle


\section{Introduction}

\begin{figure}[ht!]\anchor{typeII}\small
$$
\begin{array}{ccc}
\vcenter{\hbox{\includegraphics[height=1cm]{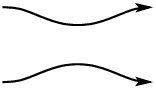}}} & \longleftrightarrow & \vcenter{\hbox{\includegraphics[height=1cm]{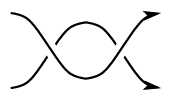}}}\\
& II_a &\\[5mm]
\vcenter{\hbox{\includegraphics[height=1cm]{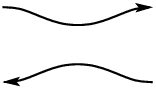}}} & \longleftrightarrow & \vcenter{\hbox{\includegraphics[height=1cm]{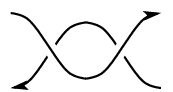}}}\\
& II_b &
\end{array}$$
\caption{The $2$ local Reidemeister moves of type $II$} \label{typeII}
\end{figure}

\begin{figure}[ht!]\anchor{typeIII}\small
$$
\begin{array}{ccccccc}
\vcenter{\hbox{\includegraphics[height=1.5cm]{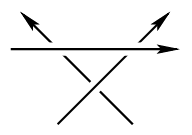}}} & \longleftrightarrow & \vcenter{\hbox{\includegraphics[height=1.5cm]{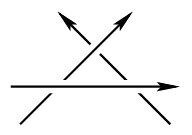}}} & 
\ \ \ \ \ \ \ & \vcenter{\hbox{\includegraphics[height=1.5cm]{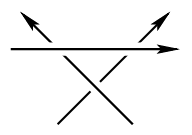}}} & \longleftrightarrow & \vcenter{\hbox{\includegraphics[height=1.5cm]{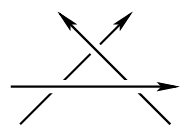}}} \\
& III_a & & & & III_b & \\[5mm]
\vcenter{\hbox{\includegraphics[height=1.5cm]{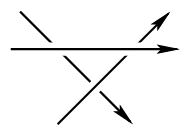}}} & \longleftrightarrow & \vcenter{\hbox{\includegraphics[height=1.5cm]{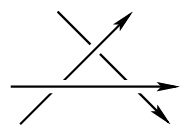}}} & 
\ \ \ \ \ \ \ & \vcenter{\hbox{\includegraphics[height=1.5cm]{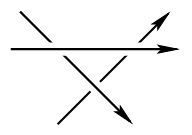}}} & \longleftrightarrow & \vcenter{\hbox{\includegraphics[height=1.5cm]{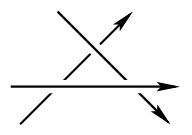}}} \\
& III_c & & & & III_d & \\[5mm]
\vcenter{\hbox{\includegraphics[height=1.5cm]{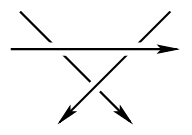}}} & \longleftrightarrow & \vcenter{\hbox{\includegraphics[height=1.5cm]{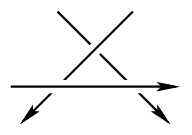}}} & 
\ \ \ \ \ \ \ & \vcenter{\hbox{\includegraphics[height=1.5cm]{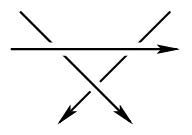}}} & \longleftrightarrow & \vcenter{\hbox{\includegraphics[height=1.5cm]{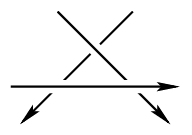}}} \\
& III_e & & & & III_f &\\[5mm]
\vcenter{\hbox{\includegraphics[height=1.5cm]{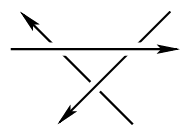}}} & \longleftrightarrow & \vcenter{\hbox{\includegraphics[height=1.5cm]{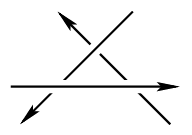}}} & 
\ \ \ \ \ \ \ & \vcenter{\hbox{\includegraphics[height=1.5cm]{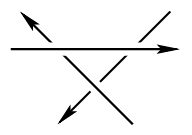}}} & \longleftrightarrow & \vcenter{\hbox{\includegraphics[height=1.5cm]{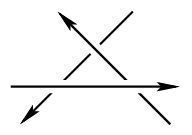}}} \\
& III_g & & & & III_h &
\end{array}$$
\caption{The $8$ local Reidemeister moves of type $III$} \label{typeIII}
\end{figure}

Links in $3$--space are usually given by diagrams and isotopies of links by sequences of Reidemeister moves (see e.g.\ \cite{Burde}). For relatively oriented links (i.e.\ up to global orientation reversing), there are exactly two different local (i.e.\ without regarding the rest of the diagram) Reidemeister moves of type $II$, shown in \figref{typeII},  and eight local Reidemeister moves of type $III$, shown in \figref{typeIII}.

\begin{defi}
The moves $II_a$, $III_a$, $III_b$, $III_c$, $III_d$, $III_e$ and $III_f$ are called braid-like Reidemeister moves. Two oriented links are called braid-like isotopic if their diagrams are related by a sequence of only braid-like Reidemeister moves.
\end{defi}

\begin{remarque}
A well-known theorem of Artin (see e.g.\ \cite{Morton}) says that two closed braids in the solid torus are isotopic if and only if they are braid-like isotopic.
\end{remarque}

It is a natural and interesting question to know in what extend certain types of Reidemeister moves can be replaced by certain other types of Reidemeister moves in an isotopy of oriented links.

Using results from chapter 1.2 in \cite{Fiedler2}, one easily proves the following lemmata.

Let $I$ denote any invariant of oriented diagrams.

\begin{lemme}
If $I$ is invariant under all Reidemeister moves of type $II$ and under just one Reidemeister move of type $III$, then it is invariant under all other Reidemeister moves of type $III$.
\end{lemme}

We have therefore to exclude some types of Reidemeister $II$ moves in order to obtain a distinction between the different types of Reidemeister $III$ moves.

\begin{lemme}\label{jeudemove}
If $I$ is invariant (only) under the braid-like Reidemeister $II$ move and under just one braid-like Reidemeister $III$ move, then it is invariant under all other braid-like Reidemeister $III$ moves.
\end{lemme}

\begin{question}
Let us consider an isotopy which uses only braid-like Reidemeister moves of type $II$ and moves of type $III_g$ and $III_h$ (the non braid-like moves of type $III$). Can the isotopy be replaced by another one which uses only braid-like Reidemeister moves of type $II$ and only exactly one of the types $III_g$ or $III_h$?
\end{question}

In other words, do the moves $III_g$ and $III_h$ become independant?

In this paper, we study oriented diagrams up to braid-like isotopies. Our main result says that the Jones polynomial can be refined in such a way that it becomes invariant in general only under braid-like isotopies. In the case of closed braids, our invariant coincides with Hoste's and Przytycki's refinement of the Jones polynomial \cite{Hoste}.

For a lightened version of our invariant, we define a corresponding Khovanov link homology. Again, in the case of closed braids, it coincides with the Khovanov homology for links in the solid torus defined by Asaeda, Przytycki and Sikora \cite{Asaeda}.

In a forthcoming paper, we construct another generalization of the
Jones polynomial for which the corresponding link homology turns out
to be invariant under move $III_h$ (as well as under the braid-like
Reidemeister $II$ move) but at the same time we have no proof (nor
counterexample) for the invariance under move $III_g$. This difficulty
sheds new light on braid-like isotopies, which form a rather natural
class of restricted isotopies of oriented links in $3$--space.

\begin{question}
Is the natural analogue of Markov's theorem still true, i.e.\ do braid-like isotopies and Markov moves (i.e.\ the four types of Reidemeister I moves in the plane) generate isotopies of oriented links in $3$--space?
\end{question}

It is a well-known open problem to know whether the Jones polynomial
detects the unknot. One can ask whether the refined Jones polynomial
or the corresponding link homology characterises the unknot. This might
be an easier question than the original one, because there are
infinitely many refined polynomials which are realised by the unknot
but we have no characterisation of these polynomials.

\section{Main results}
Our starting point is the Kauffman bracket for a framed oriented link $L$ in $3$--space. We use Kauffman's notations and terminology (see \cite{Kauffman}).

Each circle $s\co S^1 \hookrightarrow \R^2$ in a Kauffman state of a
diagram $D$ of $L$ has a piecewise orientation induced by that of
$D$. The points on $s$ where the orientation changes are called
\textit{break points}. Obviously, the number of break points on $s$ is
always even.

\begin{defi}
A circle in a Kauffman state is of type $d$ if one half of the number of break points on $s$ is odd and it is of type $h$ otherwise.
\end{defi}

Let $C$ denote the set of all configurations (up to ambiant isotopies) of unoriented embedded circles in the plane, including the empty set.

\begin{defi}
For any state $s$ of a diagram $D$, we define
$$
\begin{array}{l}
\sigma(s) = \sharp (A\textrm{--}smoothings) - \sharp (A^{-1}\textrm{--}smoothings),\\
d(s) = \sharp(d\textrm{--}circles),\\
h(s) = \sharp(h\textrm{--}circles),\\
\end{array}
$$
and $c(s) \in C$ which is the configuration of only the $h$--circles.
\end{defi}

Let $\Gamma$ denote the $\Z[A,A^{-1}]$--module over $C$.

\begin{defi}
For any diagram $D$, the bracket $\langle D \rangle_{br} \in \Gamma$ is defined by
$$
\langle D \rangle_{br}=\sum_{s\textrm{ state of }D} A^{\sigma(s)}(-A^2 - A^{-2})^{d(s)} c(s).
$$
\end{defi}

From this definition, it follows immediately:
\begin{prop}\label{premprop}
Let $D$ be a diagram and $v$ a crossing of $D$. Let $D_0$ and $D_1$ be the diagrams obtained from $D$ by smoothing $v$ respectively in the $A$--fashion and in the $A^{-1}$ one and adding break points in the natural way. Then, we have
$$
\langle D \rangle = A \langle D_0 \rangle + A^{-1} \langle D_1 \rangle.
$$
$$
\langle D \ \vcenter{\hbox{\includegraphics[height=.45cm]{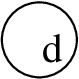}}} \rangle = (-A^2 - A^{-2}) \langle D \rangle.\leqno{\hbox{Moreover, we have}\hspace{-2cm}}
$$
\end{prop}

Let $Sei(D)$ denote the Seifert state of a diagram $D$ i.e.\ the state obtained by smoothing all the crossings with respect to the orientation of $D$, and let $w(D)$ denote the writhe of that diagram.

\begin{theo}\label{Th1}
The bracket $\langle \ . \ \rangle_{br}$ is invariant under global orientation reversing and braid-like isotopies. It has the following properties:
\begin{enumerate}
\item $(-A)^{-3w(\ .\ )}\langle \ .\ \rangle_{br}$ satisfies the skein relation of the Jones polynomial (after replacing as usual $A$ by $t^{\frac{1}{4}}$) but with infinitely many initial conditions;

\item in the case of closed braids, $(-A)^{-3w(\ .\ )}\langle \ .\ \rangle_{br}$ coincides with the invariant of Hoste and Przytycki \cite{Hoste};

\item for any diagram $D$, $\langle D \rangle_{br} =A^{w(D)} Sei(D) \ +$ terms with fewer circles in their configurations. 
\end{enumerate}
\end{theo}

\begin{remarque}
Surprisingly, we have no proof that the properties $1.$, $2.$ and $3.$ are sufficient to determine $\langle \ .\ \rangle_{br}$. The problem is that we have no classification of braid-like isotopy classes of the trivial knot (compare with Theorem 2 in \cite{Fiedler}).
\end{remarque}

\begin{question}
Is it true that for each diagram $D$ of a trivial knot, the oriented configuration of the Seifert state together with the writhe $w(D)$ determine $\langle D \rangle_{br}$?
\end{question}

\begin{remarque}
Replacing $c$ by $\chi^{|c|}$ defines a map
$$
\chi\co \Gamma \longrightarrow \Z[A,A^{-1},\chi].
$$
\end{remarque}

Obviously, for any link $L$ up to braid-like isotopy, $\langle L \rangle_{br}$ contains more information than $\chi\big(\langle L \rangle_{br}\big)$ but it turns out that the latter can be easily categorified contrary to the former.

\begin{theo}\label{Th2}
For any link $L$, there is a graded chain complex with integer coefficients $(\C_{i,j,k},d)$ with differential $d$ of tridegree $(-1,0,0)$ such that its homology groups are invariant under braid-like isotopy of the oriented link $L$.
Moreover, the bigraded Euler characteristic of this homology coincides with $(-A)^{-3w(L)}\chi\big(\langle L \rangle_{br}\big)$, that is
$$
(-A)^{-3w(L)}\chi\big(\langle L \rangle_{br}\big) = \sum_{i,j,k} (-1)^i(-A^2)^j\chi^k dim_\Q \big(\H_{i,j,k}(L)\otimes_\Z \Q\big).
$$
In the case of closed braids, the groups $\H_{i,j,k}(L)$ coincide with those defined in \cite{Asaeda}.
\end{theo}

\begin{question}
Is it possible to categorify the full bracket $\langle \ .\ \rangle_{br}$?
\end{question}

\section{Proof of theorem \ref{Th1}}

It is clear that reversing the orientation does not change break points nor $A$ and $A^{-1}$--smoothings. Thus, $\langle \ .\ \rangle_{br}$ is invariant under global orientation reversing.

From now on, we will denote break points by dots on the diagram and we will not always write brackets.

\begin{lemme}\label{lem1}
The bracket $\langle \ .\ \rangle_{br}$ is invariant under braid-like Reidemeister move of type $II$.
\end{lemme}

\begin{proof}
Calculating the Kauffman bracket, we obtain
$$
\vcenter{\hbox{\includegraphics[height=.7cm]{SchIIa2}}} \ = \ A^2 \vcenter{\hbox{\includegraphics[height=.7cm]{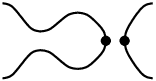}}} \ + \ \vcenter{\hbox{\includegraphics[height=.7cm]{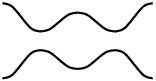}}} \ + \ \vcenter{\hbox{\includegraphics[height=.7cm]{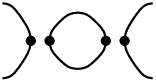}}}\ + \  A^{-2} \vcenter{\hbox{\includegraphics[height=.7cm]{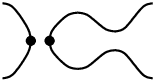}}} .
$$
The circle is of type $d$ and the usual identification $d = -A^2 - A^{-2}$ implies invariance. 
\end{proof}

\begin{remarque}
We cannot achieve invariance under the other Reidemeister $II$ move and distinguish the types $d$ and $h$ of the circles. Indeed, the calculation gives
$$
\vcenter{\hbox{\includegraphics[height=.7cm]{SchIIc2}}} \ = \ A^2 \vcenter{\hbox{\includegraphics[height=.7cm]{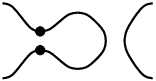}}} \ + \ \vcenter{\hbox{\includegraphics[height=.7cm]{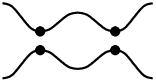}}}  \ + \ \vcenter{\hbox{\includegraphics[height=.7cm]{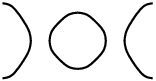}}}\ + \ A^{-2} \vcenter{\hbox{\includegraphics[height=.7cm]{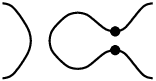}}},
$$
and the diagrams are completely different with respect to the types $d$ or $h$ of the circles!
\end{remarque}

\begin{lemme}
The bracket $\langle \ .\ \rangle_{br}$ is invariant under braid-like Reidemeister moves of type $III$.
\end{lemme}

\begin{proof}
According to lemma \ref{jeudemove}, it is sufficient to prove the invariance for the move shown in \figref{blRmIII}.

\begin{figure}[ht!]\anchor{blRmIII}
$$
\vcenter{\hbox{\includegraphics[height=1.5cm]{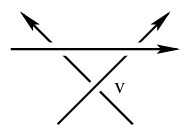}}}  \longleftrightarrow  \vcenter{\hbox{\includegraphics[height=1.5cm]{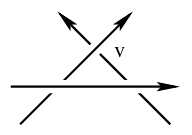}}}
$$
\caption{Braid like Reidemeister move of type $III$} \label{blRmIII}
\end{figure}

The $A$--smoothing of $v$ in the upper part of \figref{blRmIII} leads to isotopic diagrams. Hence, we have only to prove:
$$
\vcenter{\hbox{\includegraphics[height=1.05cm]{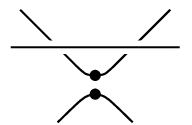}}}\ \ \ \ \ \ \ = \ \ \ \ \ \ \  \vcenter{\hbox{\includegraphics[height=1.05cm]{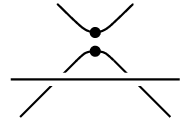}}}
$$
The calculation on the left-hand side gives
$$
\begin{array}{rcl}
\vcenter{\hbox{\includegraphics[height=1.05cm]{fig3a}}} & = &
A^2 \vcenter{\hbox{\includegraphics[height=1.05cm]{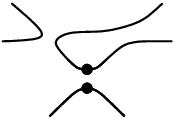}}} \ + \ \vcenter{\hbox{\includegraphics[height=1.05cm]{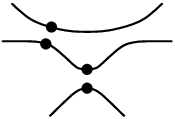}}} \ + \ \vcenter{\hbox{\includegraphics[height=1.05cm]{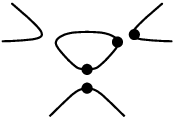}}} \ + \ A^{-2} \vcenter{\hbox{\includegraphics[height=1.05cm]{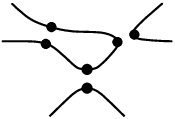}}}\\[5mm]
& = & \vcenter{\hbox{\includegraphics[height=1.05cm]{cac1}}},
\end{array}
$$
and on the right-hand side
$$
\begin{array}{rcl}
\vcenter{\hbox{\includegraphics[height=1.05cm]{fig3b}}} & = &
A^2 \vcenter{\hbox{\includegraphics[height=1.05cm]{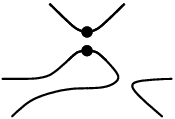}}} \ + \ \vcenter{\hbox{\includegraphics[height=1.05cm]{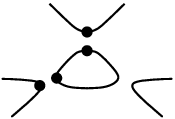}}} \ + \ \vcenter{\hbox{\includegraphics[height=1.05cm]{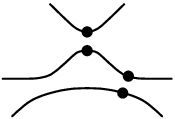}}} \ + \ A^{-2} \vcenter{\hbox{\includegraphics[height=1.05cm]{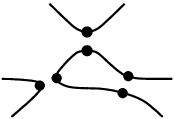}}}\\[5mm]
& = & \vcenter{\hbox{\includegraphics[height=1.05cm]{cbc3}}}.
\end{array}
$$

\begin{remarque}
Note that we are forced to distinguish the $d$ and the $h$ types by $\frac{1}{2}\sharp(\textrm{break points}) \ mod \ 2$ and we would not be able to do better than that.
\end{remarque}

Property (1) follows immediatly from proposition \ref{premprop}. There are infinitely many initial conditions because there are infinitely many diagrams of the trivial knot with the same writhe and the same Whitney index (i.e.\ that are regularly isotopic) but which are pairwise not braid-like isotopic (see \figref{cexemple}).

\begin{figure}[ht!]\anchor{cexemple}
$$\vcenter{\hbox{\includegraphics[height=1.05cm]{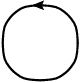}}}\ , \ \  \vcenter{\hbox{\includegraphics[height=1.05cm]{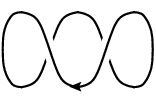}}}\ , \ \  \vcenter{\hbox{\includegraphics[height=1.05cm]{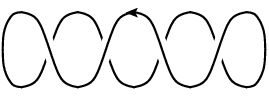}}}\ , \ \ \dots$$
\caption{Regularly, but not braid-like, isotopic trivial knots} \label{cexemple}
\end{figure}

Property (2) follows from the following simple observation: a circle in a state for a closed braid in a solid torus is contractible in the solid torus if and only if it is of type $d$. Indeed, the critical points of the restriction of the radius function on a circle correspond exactly to the break points on that circle. Evidently, the number of such critical points is congruent to $2\ mod \ 4$ if and only if the circle is contractible.

Property (3) follows from the fact that the Seifert state of a diagram $D$ contributes $A^{w(D)} Sei(D)$ to $\langle D \rangle_{br}$. Let us prove that any other state $s$ has strictly fewer circles of type $h$ than the Seifert state.

Suppose that $s$ maximizes the number of $h$--circles without being the Seifert state i.e.\ containing some break points. Then, remove as many break points as possible without changing the number of $h$--circles using the moves $1$--$4$ in figure \ref{d-break}.

\begin{figure}[ht!]\anchor{d-break}\small
$$
\begin{array}{ccccccc}
\vcenter{\hbox{\includegraphics[height=1cm]{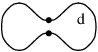}}} & \longrightarrow & \vcenter{\hbox{\includegraphics[height=1cm]{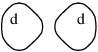}}} & \ \ \ \ \ \ \ \ \ & \vcenter{\hbox{\includegraphics[height=1cm]{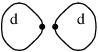}}} & \longrightarrow & \vcenter{\hbox{\includegraphics[height=1cm]{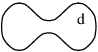}}}\\
&\textrm{move 1}&&&&\textrm{move 2}&\\[5mm]
\vcenter{\hbox{\includegraphics[height=1cm]{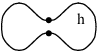}}} & \longrightarrow & \vcenter{\hbox{\includegraphics[height=1cm]{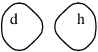}}} & \ \ \ \ \ \ \ \ \ & \vcenter{\hbox{\includegraphics[height=1cm]{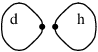}}} & \longrightarrow & \vcenter{\hbox{\includegraphics[height=1cm]{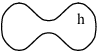}}}\\
&\textrm{move 3}&&&&\textrm{move 4}&\\[5mm]
\vcenter{\hbox{\includegraphics[height=1cm]{rbp11}}} & \longrightarrow & \vcenter{\hbox{\includegraphics[height=1cm]{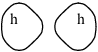}}} & \ \ \ \ \ \ \ \ \ & \vcenter{\hbox{\includegraphics[height=1cm]{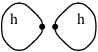}}} & \longrightarrow & \vcenter{\hbox{\includegraphics[height=1cm]{rbp22}}}\\
&\textrm{move 5}&&&&\textrm{move 6}&
\end{array}$$
\caption{Moves removing break points} \label{d-break}
\end{figure}

If a $d$--circle remains, then necessarily a move $5$ can be performed on it. As the move $5$ increases the number of $h$--circles, the resulting state must be the Seifert state by hypothesis. The result is then proven.

Otherwise, we get a configuration of only $h$--circles that we denote $s'$ with the following properties:

\renewcommand{\theenumi}{\roman{enumi}}
\renewcommand{\labelenumi}{\rm(\theenumi)}

\begin{enumerate}
\item no circle of $s'$ contains any internal breaking-smoothed crossings i.e.\ any pair of break points linked by a crossing are on distincts circles;

\item it contains some break points because, as the last move to get the Seifert state must be a move $5$ in figure \ref{d-break}, $s'$ cannot be the Seifert state by construction.
\end{enumerate}

\renewcommand{\theenumi}{\arabic{enumi}}
\renewcommand{\labelenumi}{\theenumi}

Then we construct a graph $G$ by associating a vertex to each $h$--circle of $s'$ and an edge to each pair of break points linked by a crossing. Property (i) assures that $G$ is planar and property (ii) that it contains some connected component not reduced to a point. Let's choose such a component $G'$ with $k \geq 2$ vertices and $n$ edges corresponding to a sub-state $s''$ of $s'$ with $k$ $h$--circle and $2n$ break points.

Now change the smoothing of $(k-1)$ breaking-smoothed crossing of $s''$ in such a way that the number of circle decreases each time by one. This can be done because $G'$ is connected and the removal of two linked break points corresponds to the retraction of the corresponding edge in $G'$.

Then one circle remains with $2(n-k+1)$ break points corresponding to a graph with a single vertex and $(n-k+1)$ loops. As this graph is still planar, we can find two adjacent break points linked by a crossing. Removing this two break points will create an $h$--circle with no break points. By repeating this operation, we will create one such $h$--circle for each pair of linked break points, except for the last one which will create two such $h$--circle.

At last, we get a state with $(n-k+2)$ $h$--circles. Now, each initial $h$--circles in $s''$ was containing at least one break point, therefore at least four break points in order to be of type $h$. Thus, we have
$$
2n \geq 4k,
$$
$$
n-k+2 \geq k+2.\leqno{\hbox{and consequently}\hspace{-2cm}}
$$
We have therefore increased the number of $h$--circles, so, by our hypothesis on $s$, the last state must be the Seifert state. Property (3) is thus proven.

We have thus finished the proof of the theorem.
\end{proof}

\section{Proof of theorem \ref{Th2}}

\subsection{Definition of the homology}

\begin{defi}
An enhanced state $S$ of a diagram $D$ is a Kauffman state $s$ of $D$ enhanced by an assignment of a plus or a minus sign to each of the circles of $s$.
\end{defi}

Let us fix a diagram $D$. Setting $\chi=-H^2 - H^{-2}$ and according to \cite{Viro}, we rewrite our lightened bracket as
\begin{equation}\label{bracket}
\chi\big((-A)^{-3w(D)}\langle D \rangle_{br}\big)=\!\!\!  \sum_{\substack{S\textrm{ enhanced}\\ \textrm{state of }D}} (-1)^{\frac{\sigma(S) - w(D)}{2}}(-A^2)^{\frac{\sigma(S)-3w(D)+2\tau_d(S)}{2}}(-H^2)^{\tau_h(S)},
\end{equation}
where $\tau_d$ (resp.\ $\tau_h$) is the difference between the number of pluses and minuses assigned to the $d$--circles (resp.\ $h$--circles).

\begin{defi}
For any enhanced state $S$ of $D$, we define
$$
\begin{array}{l}
i(S) = \frac{\sigma(S) - w(D)}{2}\\[2mm]
j(S) = \frac{\sigma(S)-3w(D)+2\tau_d(S)}{2}\\[2mm]
k(S) = \tau_h(S).
\end{array}
$$ 
Then, we can define $\C_{i,j,k}$ to be the $\Z$--module spanned by
$$
\{S\textrm{ enhanced state of }D\ |\ i(S) = i, j(S)=j, k(S)=k\}.
$$
\end{defi}
Strictly speaking, in accordance with the first definitions, we will not need the writhe as it is invariant under braid-like Reidemeister moves, but we will continue to use it, on the one hand for historical reasons and on the other for $i$, $j$ and $k$ to be integers.

Now, we assign an ordering of the crossings of $D$.
\begin{defi}
Let $S$ and $S'$ be two enhanced states of $D$ and $v$ a crossing of $D$. We define an incidence number $[S:S']_v$ by

\begin{itemize}
\item $[S:S']_v=1$ if the following four conditions are satisfied:
\begin{enumerate}
\item $v$ is $A$--smoothed in $S$ but $A^{-1}$--smoothed in $S'$;
\item all the other crossings are smoothed the same way in $S$ and $S'$;
\item common circles of $S$ and $S'$ are labelled the same way;
\item $j(S) = j(S')$ and $k(S) = k(S')$;
\end{enumerate}
\item otherwise, $[S:S']_v=0$.
\end{itemize}
\end{defi}

\begin{defi}
We define a differential of tridegree $(-1,0,0)$ on $(\C_{i,j,k})_{i,j,k\in \Z}$ by
$$
d(S) =  \sum_{\substack{v\textrm{ crossing}\\A\textrm{--smoothed in }S}} d_v(S)
$$
where the partial differential $d_v$ is defined by
$$
d_v(S) = (-1)^{t_{v,S}} \sum_{S'\textrm{ state of }D} [S,S']_v S',
$$
with $t_{v,S}$ the number of $A^{-1}$--smoothed crossings in $S$ labelled with numbers greater than the label of $v$.
\end{defi}

\begin{figure}[ht!]\anchor{diff}\small
\psfrag {+d}{\normalsize$\scriptstyle+d$}
\psfrag {-d}{\normalsize$\scriptstyle-d$}
\psfrag {+h}{\normalsize$\scriptstyle+h$}
\psfrag {-h}{\normalsize$\scriptstyle-h$}
$$\begin{array}{ccc}
\begin{array}{ccc}
\vcenter{\hbox{\includegraphics[height=1cm]{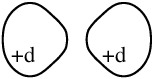}}}& \longrightarrow &\hbox{zero}\\[5mm]
\vcenter{\hbox{\includegraphics[height=1cm]{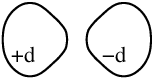}}}& \longrightarrow  &\vcenter{\hbox{\includegraphics[height=1cm]{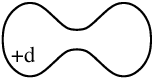}}}\\[5mm]
\vcenter{\hbox{\includegraphics[height=1cm]{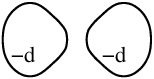}}}& \longrightarrow  &\vcenter{\hbox{\includegraphics[height=1cm]{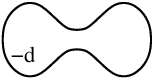}}}\\[5mm]
\vcenter{\hbox{\includegraphics[height=1cm]{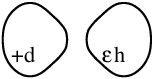}}}&  \longrightarrow  &\hbox{zero}\\[5mm]
\end{array}& \ \ &
\begin{array}{ccc}
\vcenter{\hbox{\includegraphics[height=1cm]{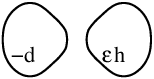}}}& \longrightarrow  &\vcenter{\hbox{\includegraphics[height=1cm]{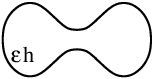}}}\\[5mm]
\vcenter{\hbox{\includegraphics[height=1cm]{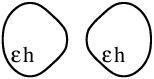}}}& \longrightarrow &\hbox{zero}\\[5mm]
\vcenter{\hbox{\includegraphics[height=1cm]{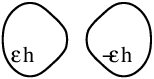}}}& \longrightarrow  &\vcenter{\hbox{\includegraphics[height=1cm]{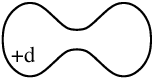}}}\\[5mm]
\end{array}
\end{array}$$\vspace{-3mm}
$$\begin{array}{ccl}
\vcenter{\hbox{\includegraphics[height=1cm]{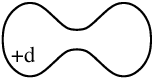}}}& \longrightarrow  &
\left\{\begin{array}{cl}
\vcenter{\hbox{\includegraphics[height=1cm]{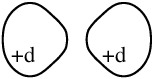}}}&*\\[5mm]
\hbox{zero}&
\end{array}\right.
\\[10mm]
\vcenter{\hbox{\includegraphics[height=1cm]{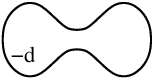}}}& \longrightarrow &
\left\{ \begin{array}{cl}
\vcenter{\hbox{\includegraphics[height=1cm]{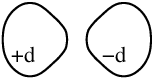}}}+\vcenter{\hbox{\includegraphics[height=1cm]{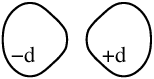}}}&*\\[5mm]
\vcenter{\hbox{\includegraphics[height=1cm]{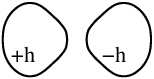}}}+\vcenter{\hbox{\includegraphics[height=1cm]{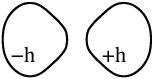}}}&\\
\end{array} \right.
\\[12mm]
\vcenter{\hbox{\includegraphics[height=1cm]{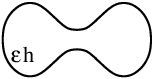}}}& \longrightarrow  &\vcenter{\hbox{\includegraphics[height=1cm]{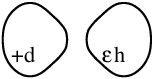}}}\\
\end{array}$$
* depending on the type of the resulting circles
\caption{Definition of the partial differential $d_v$} \label{diff}
\end{figure}

Note that $d(S)$ is the alternating sum of enhanced states obtained by switching one $A$--smoothed crossing of $S$ into a $A^{-1}$--smoothing and by locally relabeling in such a way that $j$ and $k$ are preserved and $i$ is decreased by $1$. Moreover the merge of two $d$--circles or two $h$--circles always gives a $d$--circle and the merge of a $d$--circle with an $h$--circle always an $h$--circle. Thus, one can explicit the action of $d_v$, as done in \figref{diff}.

\begin{lemme}\label{identification}
Merging rules between $d$ and $h$--circles are identical to the merging rules for the annulus in \cite{Asaeda}. 
\end{lemme}
\begin{proof}
It follows from the identification between the two definitions of $d$ and $h$--circles as given in the proof of Property (2) of theorem \ref{Th1}. 
\end{proof}

\begin{lemme}
Every two distinct partial differentials anti-commute. Thus $d^2=0$, that is $d$ is a differential. 
\end{lemme}
\begin{proof}
According to lemma \ref{identification}, one only has to check cases from $1A$--$1D$ to $5A$--$5D$ like in \cite{Asaeda} as done in the proof of lemma $4.4$ of the same article. 
\end{proof}

\begin{lemme}
The homology $\H_{br}$ associated to $d$ does not depend on the ordering of the crossings.
\end{lemme}
\begin{proof}
This follows from the proof of theorem 6.1 in \cite{Asaeda}.
\end{proof}

\begin{prop}
The extended Euler characteristic of $\H_{br}$ associated to $d$ coincides with $(-A)^{-3w(D)}\chi\big(\langle D \rangle_{br}\big)$.
\end{prop}
\begin{proof}
The fact that the extended Euler characteristic of $(\C_{i,j,k})_{i,j,k\in \Z}$ coincides with $(-A)^{-3w(D)}\chi\big(\langle D \rangle_{br}\big)$ is clear from equation \ref{bracket} and the definition of the chain complex. Moreover, passing to homology does not change the extended Euler characteristic. 
\end{proof}

Althought we still have a Viro-like exact sequence, the proof of invariance under Reidemeister moves of type $II$ and $III$ given by Asaeda, Przytycki and Sikora does not work in our context since we have no invariance under Reidemeister move of type I. Nevertheless, the original proof by Khovanov still works. We will check this using the point of view and the (intuitive) notations developed by the first author in \cite{Audoux}.

For this, we introduce the operator $.\{n\}$ (resp.\ $.[n]$) which is the global uplifting of the grading (resp.\ homological) degree by $n$.

\begin{defi}
Let $f \co \C \longrightarrow \C'$ be a morphism of chain complexes. Then the cone of $f$ is the chain complex $\C''$ defined by $\C'' = \C \oplus \C'[-1]$ and $d_{\C''}=(d_\C + f)\oplus (-d_{\C'})$.
\end{defi}

\begin{prop}{\rm\cite{Audoux}}\label{util1}\qua
Let $D$ be a diagram of a link $L$ and $v$ a crossing of $D$. Let $D_0$ and $D_1$ be the diagrams obtained from $D$ by smoothing $v$ respectively in the $A$--fashion and in the $A^{-1}$ one. Then $d_v$ defines a grading preserving map of chain complexes from $\C(D_0)$ to $\C(D_1)\{-1\}$. Moreover, $\C(D)$ is the cone of this map.
\end{prop}

\begin{prop}\label{util2}
The cone of an isomorphism between two chain complexes is acyclic.
\end{prop}

\subsection{Invariance under braid-like Reidemeister move of type $II$}

Consider two diagrams which differ only locally by a braid-like Reidemeister move of type $II$.
$$
\begin{array}{ccc}
D & \ \ \ \ \ \ \ & D'\\
\includegraphics[height=1cm]{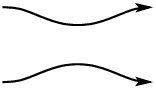} & \ \ \ \ \ \ \ & \includegraphics[height=1cm]{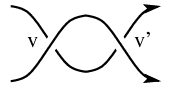}
\end{array}
$$
We order the crossings of $D$ and $D'$ in the same fashion, letting $v$ and $v'$ be the last two ones of $D'$. So we have
$$
\begin{array}{ccc}
D'_{00} &  \ \ \ \ \ \ \ & D'_{10}\\
\vcenter{\hbox{\includegraphics[height=1cm]{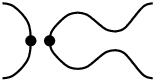}}} &  \ \ \ \ \ \ \ & \vcenter{\hbox{\includegraphics[height=1cm]{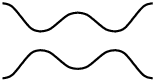}}}\\[5mm]
D'_{01} &  \ \ \ \ \ \ \ & D'_{11}\\
\vcenter{\hbox{\includegraphics[height=1cm]{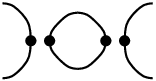}}} &  \ \ \ \ \ \ \ & \vcenter{\hbox{\includegraphics[height=1cm]{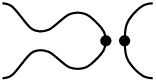}}}
\end{array}
$$
and two maps:
$$
d_{0\star} \co \vcenter{\hbox{\includegraphics[height=.7cm]{II00}}} \longrightarrow \vcenter{\hbox{\includegraphics[height=.7cm]{II01}}}
$$
$$
d_{\star 1} \co  \vcenter{\hbox{\includegraphics[height=.7cm]{II01}}}\longrightarrow \vcenter{\hbox{\includegraphics[height=.7cm]{II11}}}
$$

\begin{lemme}
The morphisms $d_{0\star}$ and $d_{\star 1}$ are respectively injective and surjective.
\end{lemme}

\begin{proof}
This can be checked by direct calculus on canonical basis. 
\end{proof}

\begin{lemme}\label{firstproof}
The cone of
$$
\xymatrix @!0 @C=3cm @R=1.5cm {
\cdots  &  Im(d^i_{\star 0}+ d^i_{0 \star})\{-1\} \ar[l]_(.68){d^{i-1}_{D'}} & \cdots \ar[l]_(.26){d^i_{D'}}\\
\cdots  & \Z\left\{\vcenter{\hbox{\includegraphics[height=.7cm]{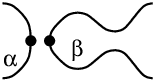}}}\right\}  \ar[u]^{d^i_{\star 0}+ d^i_{0 \star}} \ar[l]_(.6){d^{i-1}_{D'_{00}}}& \cdots \ar[l]_(.33){d^i_{D'_{00}}}\\
}
$$
is an acyclic graded subcomplex of $\C(D')$ denoted $\C_1$.
\end{lemme}
\begin{proof}
To prove this lemma, we have to check the four following points:
\begin{enumerate}
\item stability of the upper and the lower chain complexes under their differentials;
\item grading-preserving;
\item inclusion in $\C(D')$;
\item acyclicity.
\end{enumerate}
Stability is clear for the lower line. For the upper one, as $d_{0\star}$ is injective, the map $(d^i_{\star 0}+ d^i_{0 \star})$ has a right inverse on its image, and by anti-commutativity of distinct partial differentials we have for any $S$ in $Im(d^i_{\star 0}+ d^i_{0 \star})\{-1\}$
$$
d^{i-1}_{D'}(S) = -(d^{i-1}_{\star 0}+ d^{i-1}_{0 \star})\circ d^{i-1}_{D'_{00}} \circ (d^i_{\star 0}+ d^i_{0 \star})^{-1}(S) \in Im(d^{i-1}_{\star 0}+ d^{i-1}_{0 \star})\{-1\}.
$$
The second point is clear by definition.

Inclusion of spaces is obvious and the identification of the differentials follows from proposition \ref{util1}.

The last point is a consequence of proposition \ref{util2}. 
\end{proof}

\begin{lemme}
The cone of
$$
\xymatrix @!0 @C=3cm @R=1.5cm {
\cdots &  \Z\left\{\vcenter{\hbox{\includegraphics[height=.7cm]{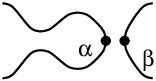}}}\right\} \{-1\} \ar[l]_(.65){d^{i-1}_{D'}}  & \cdots \ar[l]_(.27){d^i_{D'}}\\
\cdots   & \Z\left\{\vcenter{\hbox{\includegraphics[height=.7cm]{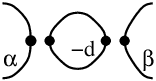}}}\right\}  \ar[u]^{d^i_{\star 1}} \ar[l]_(.61){d^{i-1}_{D'_{01}}} & \cdots \ar[l]_(.35){d^i_{D'_{01}}}
}
$$
is an acyclic graded subcomplex of $\C(D')$ denoted $\C_2$.
\end{lemme}

\begin{proof}
The stability of the upper line is clear. The differential $d^i_{D'_{01}}$ do not change the smoothing of the only two crossings contiguous to the middle circle labeled by $-d$, so this one is let unchanged and the lower line is stable too.

Arguments of lemma \ref{firstproof} are still valid to prove the last three points. 
\end{proof}

\begin{defi}
We define the map:
$$
\begin{array}{rccc}
.\otimes v_- \co & \vcenter{\hbox{\includegraphics[height=.85cm]{II11}}} & \longrightarrow & \vcenter{\hbox{\includegraphics[height=.85cm]{II01}}}\\[4mm]
& \vcenter{\hbox{\includegraphics[height=.7cm]{Inj2}}} & \mapsto & -\vcenter{\hbox{\includegraphics[height=.7cm]{Inj1}}}
\end{array}
$$
\end{defi}

\begin{remarque} $\ $
The map $.\otimes v_-$ is grading-preserving and is a right inverse for $d_{\star 1}$.
\end{remarque}

\begin{lemme}
The complex
$$
\xymatrix @!0 @C=6cm @R=2cm {
\cdots  &  \Z\left\{\vcenter{\hbox{\includegraphics[height=.7cm]{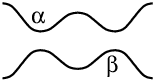}}}-\big( d_{1 \star}(\vcenter{\hbox{\includegraphics[height=.7cm]{isomII}}}) \big)\otimes v_- \right\} \ar[l]_(.7){d^{i-1}_{D'}} & \cdots \ar[l]_(.23){d^i_{D'}}}
$$
is a graded subcomplex of $\C(D')$ denoted $\C_3$, isomorphic to $\C(D)$ via $\psi_{II}$ defined by:
$$
\psi_{II}\left(\vcenter{\hbox{\includegraphics[height=.7cm]{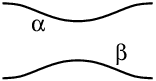}}} \right)=\vcenter{\hbox{\includegraphics[height=.7cm]{isomII}}}-\big( d_{1 \star}(\vcenter{\hbox{\includegraphics[height=.7cm]{isomII}}}) \big)\otimes v_-
$$
\end{lemme}
\begin{proof}
To begin, we split the differential into the sum of the underlying partial differentials. First consider $d_c$ where $c \neq v,v'$. As $d_c$ let unchanged the middle circle in $\vcenter{\hbox{\includegraphics[height=.5cm]{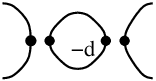}}}$ and as $.\otimes v_-$ remove one $A^{-1}$--smoothed crossing labeled with a number greater than the label of $c$, we have:
$$
\begin{array}{rcl}
d_c \Big( \big( d_{1 \star}(\vcenter{\hbox{\includegraphics[height=.7cm]{isomII}}}) \big)\otimes v_- \Big) & = & -\Big( d_c \big( d_{1 \star}(\vcenter{\hbox{\includegraphics[height=.7cm]{isomII}}}) \big)\Big) \otimes v_-\\[5mm]
& = & \Big( d_{1 \star} \big( d_c (\vcenter{\hbox{\includegraphics[height=.7cm]{isomII}}}) \big) \Big) \otimes v_-
\end{array}
$$
So we have:
$$
d_c \Big( \vcenter{\hbox{\includegraphics[height=.7cm]{isomII}}}-\big( d_{1 \star}(\vcenter{\hbox{\includegraphics[height=.7cm]{isomII}}}) \big)\otimes v_-  \Big) = d_c\Big(\vcenter{\hbox{\includegraphics[height=.7cm]{isomII}}}\Big)-\Big( d_{1 \star}\big(d_c(\vcenter{\hbox{\includegraphics[height=.7cm]{isomII}}})\big) \Big)\otimes v_-
$$
Then, $d_{1 \star}$ and $d_{\star 1}$ act respectively on  $\vcenter{\hbox{\includegraphics[height=.5cm]{isomII}}}$ and $\big( d_{1 \star}(\vcenter{\hbox{\includegraphics[height=.5cm]{isomII}}}) \big)\otimes v_-$. But as $.\otimes v_-$ is a right inverse for $d_{\star 1}$, we have:
$$
\begin{array}{rcl}
d_{1 \star}\Big( \vcenter{\hbox{\includegraphics[height=.7cm]{isomII}}} \Big) - d_{\star 1}\Big(  \big( d_{1 \star}(\vcenter{\hbox{\includegraphics[height=.5cm]{isomII}}}) \big)\otimes v_-\Big) & = & d_{1 \star}( \vcenter{\hbox{\includegraphics[height=.7cm]{isomII}}} ) - d_{1 \star}(\vcenter{\hbox{\includegraphics[height=.5cm]{isomII}}})\\[5mm]
& = & 0
\end{array}
$$
Not only this two points show that the subcomplex is stable, but also, as, for $c \neq v,v'$, $d_c$ acts the same on   
$\vcenter{\hbox{\includegraphics[height=.5cm]{isomII2}}}$ and $\vcenter{\hbox{\includegraphics[height=.5cm]{isomII}}}$, they show that $\psi_{II}$ is a chain complexes map. Then, it is clearly a grading-preserving isomorphism since $\vcenter{\hbox{\includegraphics[height=.5cm]{isomII}}}$ has only two more crossings than $\vcenter{\hbox{\includegraphics[height=.5cm]{isomII2}}}$, one $A$--smoothed and the other $A^{-1}$--smoothed. 
\end{proof}

\begin{lemme}\label{lastlemme}
$\C (D') \simeq \C_1 \oplus \C_2 \oplus \C_3$.
\end{lemme}

\begin{proof}
As $\C_1$, $\C_2$ and $\C_3$ are subcomplexes of $\C(D')$, differentials can already been identified. Thus, one only has to prove it at the level of modules.

The subcomplex $\C_2$ contains all elements of the form  
{\psfraga <-3.8pt,0pt> {-d}{{\small$\scriptscriptstyle {-}\!d$}}
$\vcenter{\hbox{\includegraphics[height=.5cm]{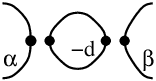}}}$,} 
thus we can clean all parasite terms in $\C_3$ and get
$$
\C_2 \oplus \C_3 \simeq \vcenter{\hbox{\includegraphics[height=.7cm]{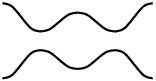}}} \oplus \vcenter{\hbox{\includegraphics[height=.7cm]{SchII01}}}\oplus \vcenter{\hbox{\includegraphics[height=.7cm]{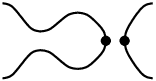}}}.
$$
Moreover,
\psfraga <-1.5pt,0pt> {h}{$\scriptscriptstyle h$}
$$
Im(d^i_{\star 0}+ d^i_{0 \star}) = \Z \left\{ \vcenter{\hbox{\includegraphics[height=.7cm]{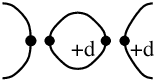}}} + \alpha_1, \vcenter{\hbox{\includegraphics[height=.7cm]{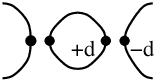}}} + \vcenter{\hbox{\includegraphics[height=.7cm]{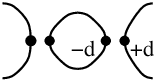}}} + \alpha_2, \vcenter{\hbox{\includegraphics[height=.7cm]{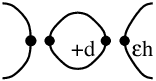}}} + \alpha_3 \right\},
$$
with $\alpha_i \in \vcenter{\hbox{\includegraphics[height=.5cm]{SchII10}}}$.
So, with $\C_2 \oplus \C_3$ we can again clean and get
$$
Im(d^i_{\star 0}+ d^i_{0 \star})\oplus \C_2 \oplus \C_3 \simeq \vcenter{\hbox{\includegraphics[height=.7cm]{SchII10}}} \oplus \vcenter{\hbox{\includegraphics[height=.7cm]{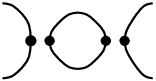}}}\oplus \vcenter{\hbox{\includegraphics[height=.7cm]{SchII11}}}.
$$
Now $\C_1 = \vcenter{\hbox{\includegraphics[height=.5cm]{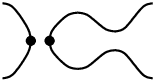}}}\oplus Im(d^i_{\star 0}+ d^i_{0 \star})$, so the lemma is proven.
\end{proof}

Finally, we achieve the invariance under star-like Reidemeister move of type $II$ by passing to homology in lemma \ref{lastlemme}.

\subsection{Invariance under braid-like Reidemeister moves of type $III$}
According to lemma \ref{jeudemove}, it is sufficient to prove the invariance for two diagrams which differ only locally by the following braid-like Reidemeister move of type $III$.
{\small$$
\begin{array}{ccc}
D & \ \ \ \ \ \ \ & D'\\
\includegraphics[height=1.5cm]{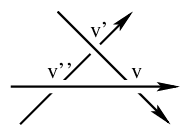} & \ \ \ \ \ \ \ & \includegraphics[height=1.5cm]{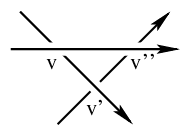}.
\end{array}
$$}%
We order the crossings of $D$ and $D'$ in the same fashion, letting the three crossings involved in the Reidemeister move be the last three ones. Thus we have:
{\small$$
\begin{array}{ccc}
\begin{array}{ccc}
D_{000}& \ \ \ \ & D_{001}\\
\vcenter{\hbox{\includegraphics[height=1.4cm]{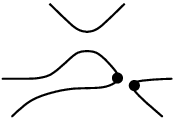}}} &  \ \ \ \ & \vcenter{\hbox{\includegraphics[height=1.4cm]{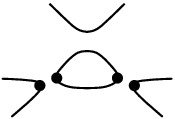}}}
\end{array}
& \ \ \ \ \ \ &
\begin{array}{ccc}
D'_{000}& \ \ \ \ & D'_{001}\\
\vcenter{\hbox{\includegraphics[height=1.4cm]{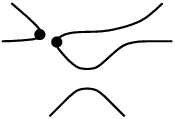}}} &  \ \ \ \ & \vcenter{\hbox{\includegraphics[height=1.4cm]{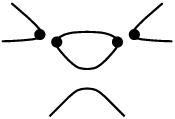}}}
\end{array}\\[10mm]
\begin{array}{ccc}
D_{010}& \ \ \ \ & D_{011}\\
\vcenter{\hbox{\includegraphics[height=1.4cm]{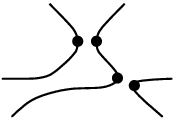}}} &  \ \ \ \ & \vcenter{\hbox{\includegraphics[height=1.4cm]{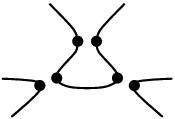}}}
\end{array}
& \ \ \ \ \ \ &
\begin{array}{ccc}
D'_{010}& \ \ \ \ & D'_{011}\\
\vcenter{\hbox{\includegraphics[height=1.4cm]{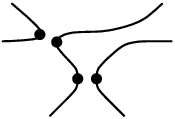}}} &  \ \ \ \ & \vcenter{\hbox{\includegraphics[height=1.4cm]{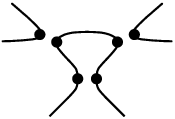}}}
\end{array}\\[10mm]
D_{1\bullet \bullet}& \ \ \ \ \ \ & D'_{1\bullet \bullet}\\
\vcenter{\hbox{\includegraphics[height=1.4cm]{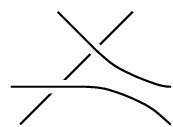}}} &  \ \ \ \ \ \ & \vcenter{\hbox{\includegraphics[height=1.4cm]{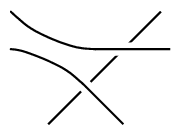}}}
\end{array}
$$}%
We can, immediatly, note that $D_{1\bullet \bullet}\simeq D'_{1\bullet \bullet}$ and $D_{010}\simeq D'_{010}$. With respect to $D$, we have two maps:
$$
d_{00\star} \co \vcenter{\hbox{\includegraphics[height=1.05cm]{III00b}}} \longrightarrow \vcenter{\hbox{\includegraphics[height=1.05cm]{III01b}}}
$$
$$
d_{0\star 1} \co  \vcenter{\hbox{\includegraphics[height=1.05cm]{III01b}}}\longrightarrow \vcenter{\hbox{\includegraphics[height=1.05cm]{cad2}}}
$$

\begin{lemme}
The maps $d_{00\star}$ and $d_{0\star 1}$ are, respectively injective and surjective.
\end{lemme}

Then, the rest follows as in the case of the Reidemeister move of type $II$.

\begin{lemme}
The cone of
$$
\xymatrix @!0 @C=4cm @R=2cm {
\cdots &  Im(d^i_{\star 00}+ d^i_{0 \star 0} + d^i_{00\star})\{-1\} \ar[l]_(.72){d^{i-1}_D} & \cdots \ar[l]_(.2){d^i_D}\\
\cdots & \Z\left\{\vcenter{\hbox{\includegraphics[height=1.05cm]{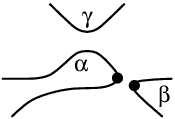}}}\right\} \ar[u]^{d^i_{\star 00}+ d^i_{0 \star 0}+d^i_{00\star}} \ar[l]_(.58){d^{i-1}_{D_{000}}} & \cdots \ar[l]_(.37){d^i_{D_{000}}}\\
}
$$
is an acyclic graded subcomplex of $\C(D)$ denoted $\C_1$.
\end{lemme}

\begin{lemme}
The cone of
$$
\xymatrix @!0 @C=4,2cm @R=2cm {
\cdots &  \Z\left\{\vcenter{\hbox{\includegraphics[height=1.05cm]{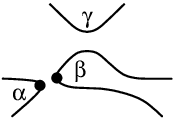}}}+\vcenter{\hbox{\includegraphics[height=1.05cm]{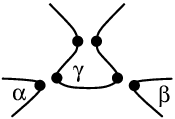}}}\right\} \{-1\} \ar[l]_(.75){d^{i-1}_D} & \cdots \ar[l]_(.17){d^i_D} \\
\cdots  & \Z\left\{\vcenter{\hbox{\includegraphics[height=1.05cm]{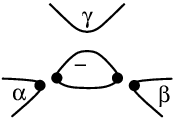}}}\right\}  \ar[u]^{d^i_{\star 01}+d^i_{0\star 1}} \ar[l]_(.57){d^{i-1}_{D_{001}}} & \cdots \ar[l]_(.37){d^i_{D_{001}}}
}
$$
is an acyclic graded subcomplex of $\C(D')$ denoted $\C_2$.
\end{lemme}

\begin{defi}
We define the map
$$
\begin{array}{rccc}
.\otimes v_- \co & \vcenter{\hbox{\includegraphics[height=1.275cm]{cad2}}} & \longrightarrow & \vcenter{\hbox{\includegraphics[height=1.275cm]{III01b}}}\\[4mm]
& \vcenter{\hbox{\includegraphics[height=1.05cm]{surjIII2}}} & \mapsto & -\vcenter{\hbox{\includegraphics[height=1.05cm]{surjIII}}}.
\end{array}
$$
\end{defi}

\begin{remarque} $\ $
The map $.\otimes v_-$ is grading-preserving and is a right inverse for $d_{0\star 1}$.
\end{remarque}

\begin{lemme}
The complex
$$
\xymatrix @!0 @C=6.2cm @R=2cm {
\cdots &  \Z\left\{\vcenter{\hbox{\includegraphics[height=1.05cm]{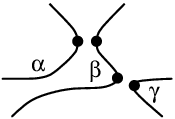}}}-\big( d_{01 \star}(\vcenter{\hbox{\includegraphics[height=1.05cm]{isomIII}}}) \big)\otimes v_- \right\}\oplus \vcenter{\hbox{\includegraphics[height=1.275cm]{passeb}}} \ar[l]_(.8){d^{i-1}_D} & \cdots \ar[l]_(.15){d^i_D} }
$$
is a graded subcomplex of $\C(D)$ denoted $\C_3$.
\end{lemme}

\begin{lemme}
$\C(D) \simeq \C_1 \oplus \C_2 \oplus \C_3$.
\end{lemme}

With the same reasoning on $D'$, we get:
\begin{lemme}
$\C(D') \simeq \C'_1 \oplus \C'_2 \oplus \C'_3$ with $\C'_1$ and $\C'_2$ acyclic and
$$
\C'_3 = \Z\left\{\vcenter{\hbox{\includegraphics[height=1.05cm]{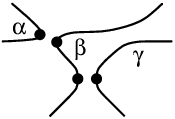}}}-\big( d_{01 \star}(\vcenter{\hbox{\includegraphics[height=1.05cm]{isomIIIb}}}) \big)\otimes v_- \right\}\oplus \vcenter{\hbox{\includegraphics[height=1.275cm]{passe}}}
$$
\end{lemme}

Finally, we achieve the invariance under braid-like Reidemeister moves
of type $III$ thanks to the following lemma.

\begin{lemme}
The chain complexes $\C_3$ and $\C'_3$ are isomorphic via $\psi_{III}$ defined on $\vcenter{\hbox{\includegraphics[height=1cm]{cad1}}}$ by
$$
\psi\left(\vcenter{\hbox{\includegraphics[height=1.05cm]{isomIII}}}-\big( d_{01 \star}(\vcenter{\hbox{\includegraphics[height=1.05cm]{isomIII}}}) \big)\otimes v_- \right)=\vcenter{\hbox{\includegraphics[height=1.05cm]{isomIIIb}}}-\big( d_{01 \star}(\vcenter{\hbox{\includegraphics[height=1.05cm]{isomIIIb}}}) \big)\otimes v_-
$$
and on $\vcenter{\hbox{\includegraphics[height=1cm]{passeb}}}$ in the obvious way.
\end{lemme}

\textbf{Acknowledgement}\qua We are grateful to the referee for numerous and useful remarks.



\Addresses\recd

\end{document}